\documentclass[12pt]{article}

\def\sqr#1#2{{\vcenter{\vbox{\hrule height.#2pt
              \hbox{\vrule width.#2pt height#1pt \kern#1pt \vrule width.#2pt}
              \hrule height.#2pt}}}}
\def\signed #1{{\unskip\nobreak\hfil\penalty50
              \hskip2em\hbox{}\nobreak\hfil#1
              \parfillskip=0pt \finalhyphendemerits=0 \par}}
\def\endpf{\signed {$\sqr69$}}

\def\dbR{{\mathop{\rm l\negthinspace R}}}
\def\dbC{{\mathop{\rm l\negthinspace\negthinspace\negthinspace C}}}

\def\3n{\negthinspace \negthinspace \negthinspace }
\def\2n{\negthinspace \negthinspace }
\def\1n{\negthinspace }

\def\dbC{{\mathop{\rm l\negthinspace\negthinspace\negthinspace C}}}

\def\dbE{{\mathop{\rm l\negthinspace E}}}
\def\dbF{{\mathop{\rm l\negthinspace F}}}

\def\ds{\displaystyle}

\def\dbN{{\mathop{\rm l\negthinspace N}}}

\def\dbR{{\mathop{\rm l\negthinspace R}}}
\def\={\buildrel \triangle \over =}

\def\a{\alpha}
\def\b{\beta}
\def\g{\gamma}
\def\d{\delta}

\def\l{\lambda}

 \def\n{\nabla}

\def\t{\times}

\def\th{\theta}

\def\u{\upsilon}
\def\i{\infty}
\def\ns{\noalign{\ss} }

\def\D{\Delta}

\def\L{\Lambda}
\def\Si{\Sigma}

\def\O{\Omega}

\def\cF{{\cal F}}

\def\cL{{\cal L}}

\def\cl{{\cal l}}
\def\no{\noindent}

\def\ss{\smallskip}
\def\ms{\medskip}
\def\bs{\bigskip}
\def\q{\quad}
\def\qq{\qquad}

\def\pa{\partial}

\def\cd{\cdot}

\def\div{\hbox{\rm div$\,$}}

\def\cl{\overline}

\def\Re{{\mathop{\rm Re}\,}}
\def\Im{{\mathop{\rm Im}\,}}

\def\|{\Big |}
\def\({\Big (}
\def\){\Big )}
\def\[{\Big[}
\def\]{\Big]}
\def\be{\begin{equation}}
\def\bel{\begin{equation}\label}
\def\ee{\end{equation}}
\def\bt{\begin{theorem}}
\def\bcd{\begin{condition}}
\def\ecd{\end{condition}}
\def\et{\end{theorem}}
\def\bc{\begin{corollary}}
\def\ec{\end{corollary}}
\def\bde{\begin{definition}}
\def\ede{\end{definition}}
\def\bl{\begin{lemma}}
\def\el{\end{lemma}}
\def\bp{\begin{proposition}}
\def\ep{\end{proposition}}
\def\br{\begin{remark}}
\def\er{\end{remark}}
\def\ba{\begin{array}}
\def\ea{\end{array}}
\def\ed{\end{document}}
\def\ns{\noalign{\ms}}
\def\ds{\displaystyle}

\def\square#1{\vbox{\hrule\hbox{\vrule height#1%
     \kern#1\vrule}\hrule}}
\def\rectangle#1#2{\vbox{\hrule\hbox{\vrule height#1%
     \kern#2\vrule}\hrule}}

\font\tenbb=msbm10 \font\sevenbb=msbm7 \font\fivebb=msbm5

\newfam\bbfam
\scriptscriptfont\bbfam=\fivebb \textfont\bbfam=\tenbb
\scriptfont\bbfam=\sevenbb

\newtheorem{lemma}{Lemma}[section]
\newtheorem{remark}{Remark}[section]

\newtheorem{theorem}{Theorem}[section]
\newtheorem{corollary}{Corollary}[section]

\newtheorem{definition}{Definition}[section]
\newtheorem{proposition}{Proposition}[section]
\newtheorem{condition}{Condition}[section]

\makeatletter
   
   \@addtoreset{equation}{section}
\makeatother

\begin{document}
\title{\bf A weighted  identity for stochastic partial differential
operators and its applications   \thanks{This work is partially
supported by the NSF of China under grants  11322110,  11371084 and 11171060.}}

\author{Xiaoyu Fu\thanks{School of Mathematics, Sichuan University, Chengdu 610064, China. E-mail address: rj$\_$xy@163.com.}  \ and \ Xu Liu\thanks{School of Mathematics and Statistics, Northeast Normal
University, Changchun 130024, China. E-mail address:
liux216@nenu.edu.cn.}}

\date{}

\maketitle

\begin{abstract}
\no  In this paper,  a pointwise weighted identity for some stochastic partial differential operators (with complex principal parts) is established. This identity presents a unified approach in studying the controllability, observability and inverse problems for some deterministic/stochastic partial differential equations. Based on this identity, one can deduce all the known Carleman estimates and observability results, for some  deterministic partial differential equations, stochastic heat equations, stochastic Schr\"odinger equations and stochastic transport equations.  Meanwhile,  as its new application,   we study an inverse problem for  linear stochastic complex Ginzburg-Landau equations.
 \end{abstract}

\bs

\no{\bf 2010 Mathematics Subject Classification}.  Primary 93B05; Secondary 93B07, 93C20
\bs

\no{\bf Key Words}.  Stochastic partial differential equations,  controllability, observability, inverse problems, global Carleman
estimate.

\ms

\section{Introduction}

\ms

This paper is devoted to a pointwise weighted identity for  a class of stochastic partial differential operators. Based on this identity,    one can derive  global Carleman estimates for deterministic/stochastic partial differential equations of different type.   This idea first came from
 the Russian literature \cite{LRS},  in which
some unique continuation results were established, based on suitable Carleman estimates.
Carleman estimates were originally introduced by T.~Carleman  in 1939.    They were   energy estimates  with
exponential weights,  and  established in \cite{C}  to prove   a strong unique continuation
property for some elliptic  equations   in
dimension two.  Up to now,   Carleman estimates have become  a powerful tool in  studying  deterministic/stochastic partial differential equations, and the  related control  and inverse problems.  For example,
this type of weighted energy estimates  was used to
study the unique continuation property of partial differential equations (\cite{Ke}),  the uniqueness and stability of Cauchy problems (\cite{BK, Ca, H}),  inverse
problems (\cite{BP, I, K}) and the controllability $(\cite{F, fur, Y, Z, zua})$, respectively.

\ms

Roughly speaking,   a Carleman estimate for the Laplacian operator is  an inequality of the form:
  $$
 |e^{\l\phi}f|_{L^2(G)}\le C|e^{\l\phi}\D f|_{L^2(G)}, \quad \mbox{for any }f\in C_0^2(G),
  $$
where $G$ is a nonempty open subset of $\dbR^n$ with a smooth boundary,  $\phi$ is a suitable weighted function, and  $C$ is a constant,  independent of the parameter $\l$ which may tend to  $+\i$.

\medskip

 In what follows, we give two simple  examples to introduce the basic idea of  establishing Carleman estimates.

\ms

\noindent{ \bf Example 1. The stability of an ordinary differential system }

\medskip

Consider the following ordinary differential system:
 \be\label{ODE}
  \left\{\ba{ll} \dot{x}(t)=a(t)x(t)\quad
 t\in (0,T), \\
 \ns\ds x(0)=x_{0}, \ea
 \right.
 \ee
 where $x_0\in \dbR^{n}$ and $a(\cdot)\in L^\infty(0, T; \dbR^{n\times n})$.
Then  for any
$\l>0$, by the first equation of (\ref{ODE}), we have that
 \bel{CD}\qq
 \frac{d}{dt}\(e^{-\l t}|x(t)|^{2}\)=-\l e^{-\l t}|x(t)|^{2}+2e^{-\l t}x(t)\cd \dot{x}(t)
 \leq [2|a(t)|_{\dbR^{n\times n}}-\l]e^{-\l t}|x(t)|^{2}.
 \ee
Choosing a sufficiently large $\l$,
one can obtain that
 $$
 |x(t)|\leq e^{\l t}|x_{0}|\leq e^{\l T}|x_0|,\quad\forall\ t  \in[0,T].
 $$
The key of this proof  for the stability of (\ref{ODE}) is   the following
identity:
 \bel{0e2}
 2e^{-\l t} x(t)\cd \dot{x}(t)=\frac{d}{dt}\(e^{-\l t}|x(t)|^{2}\)+\l e^{-\l t}|x(t)|^{2}.
 \ee
 (\ref{0e2}) can be viewed as
 a pointwise weighted identity for the principal
operator $\dot{x}(t)$ of (\ref{ODE}).  After multiplied  by a  multiplier $2e^{-\lambda t}x(t)$,   the principal operator is rewritten as a sum of a ``divergence" term
$\displaystyle\frac{d}{dt}(e^{-\l t}|x(t)|^{2})$ and an ``energy" term $\l e^{-\l
t}|x(t)|^{2}$. By choosing a sufficiently large parameter $\l$,    the  undesired lower order term $2|a(t)|_{\dbR^{n\times n}}e^{-\l t}|x(t)|^{2}$   with respect to $\lambda$ can be absorbed.

 \ms

\noindent { \bf Example 2. A  Carleman estimate for  first order  differential operators}

 \medskip

 For any fixed $\g_0\in C(\cl{G})$ and $\g\in \[C^1(\cl{G})\]^n$,   consider the following first order  differential operator:
 \bel{ccf1}
 \cL(x,D)u=\g\cdot\n u+\g_0u,\q\forall\  x\in \cl{G}.
 \ee
Set
 \bel{phi}
\phi(x)=|x-x_0|^2, \mbox{ for some }x_0\in\dbR^n.
 \ee
Then, we have the following known Carleman estimate for the operator (\ref{ccf1}).
 \bl\label{stir}   Assume that for  $x_0\in\dbR^n\setminus\cl{G}$ and a positive constant $c_0$,
  \bel{ccf2}
  \g(x)\cdot(x-x_0)\le - c_0, \q\mbox{ in } \cl{G}.
  \ee
Then  there exist  constants $\l^*>0$ and $C>0$, so that for any  $\l \ge\l^*$,
      \bel{ccf3}
\l\int_G  e^{2\l\phi}  u^2dx  \le C\int_Ge^{2\l\phi}|\mathcal{L}(x,D)u|^2dx,
   \ee
  for any  $u\in C_0^1(G)$.
 \el

 \ss

 \noindent {\bf Proof of Lemma \ref{stir}. }
For any  $\l>0$, put
$$
\ell(x)=\l \phi(x)\quad\mbox{ and }\quad \th=e^{\l\phi},
$$
where $\phi$ is given by (\ref{phi}). Then by (\ref{phi}),  it is easy to check that
 \begin{equation}\label{fst}(\theta^2 u) \gamma\cdot\nabla u=
 \theta^2 \gamma\cdot\nabla\Big(\ds\frac{1}{2}u^2\Big)=\div\Big(\ds\frac{1}{2}\theta^2 u^2\gamma\Big)
  -\theta^2\Big[\ds\frac{1}{2}\div \gamma+2\lambda\gamma\cdot(x-x_0)\Big]u^2.
 \end{equation}
This implies that there exists a constant $C_1>0$ such that
 \begin{eqnarray*}
 &&(\theta^2 u) \mathcal{L}(x, D)u=
\div\Big(\ds\frac{1}{2}\theta^2 u^2\gamma\Big)
  -\theta^2\Big[\ds\frac{1}{2}\div \gamma+2\lambda\gamma\cdot(x-x_0)-\gamma_0\Big]u^2\\
 &&\geq \div\Big(\ds\frac{1}{2}\theta^2 u^2\gamma\Big)
  -\theta^2\Big[2\lambda\gamma\cdot(x-x_0)+C_1\Big]u^2.
\end{eqnarray*}
By (\ref{ccf2}), integrating the above inequality in $G$ and  choosing  a sufficiently large $\l$,  we can get  the desired
 estimate (\ref{ccf3}). \endpf

\medskip

\medskip

 \noindent The key of this  proof of  Lemma \ref{stir}  is  the identity (\ref{fst}).
 It  can be viewed as
 a pointwise weighted identity for the principal
operator $\gamma\cdot\nabla u$ of $\mathcal{L}(x, D)$.  After multiplied  by a  multiplier $\theta^2 u$,   the principal operator is rewritten as a sum of a ``divergence" term
$\div\Big(\ds\frac{1}{2}\theta^2 u^2\gamma\Big)$ and an ``energy" term $ -\theta^2\Big[\ds\frac{1}{2}\div \gamma+2\lambda\gamma\cdot(x-x_0)\Big]u^2$. By choosing a sufficiently large parameter $\l$,    the  undesired lower order term $\ds\frac{1}{2}\div \gamma$   with respect to $\lambda$ can be absorbed.

\ms

From the above two examples,  one can find that the key of proving Carleman estimates is  to establish a suitable pointwise weighted identity for principal operators of  differential equations.  Notice that in $\cite{F}$,  a pointwise weighted identity for the following deterministic partial differential operator was established: $$Lw=(\a+i\b)w_t+\sum\limits_{j, k=1}^{n}(a^{j k}w_{x_j})_{x_k},$$
where $\a$, $\b$ and $a^{j k}(\cdot)$ $(k, j=1, 2, \cdots, n)$ are suitable real-valued functions,  and $i$ is the imaginary unit.  This identity presented  a unified approach of deducing global Carleman estimates  for many deterministic partial differential equations of different type.  A natural problem is whether one can get the counterpart for stochastic partial differential equations. As far as we know,  there exist few  works on global  Carleman estimates for stochastic partial differential equations. We refer to  \cite{Lu1, Lu2, Z} for some known results in this respect.  However,    there is not any  known Carleman estimate for  general stochastic  partial differential operators with complex principal parts.
In this paper,  we mainly present a pointwise weighted identity for the following stochastic partial differential operator:
$$
\cL w=a_0dw-(a+ib)\sum\limits_{j, k=1}^{n}(a^{j k}w_{x_j})_{x_k}dt+{\bf b_0}\cdot \n wdt,
$$
where $a_0,  a, b\in\dbR$ and ${\bf b_0}\in\dbR^n$.  The operator $\mathcal{L}$ may include some deterministic/stochastic partial differential operators of different type.
Based on a  pointwise weighted   identity for this operator,  we develop a unified approach of establishing  global Carleman estimates   for stochastic heat equations, stochastic
Schr\"{o}dinger equations, stochastic transport equations and linear stochastic complex Ginzburg-Landau equations. As  applications of this identity, one can also study some inverse problems of these different stochastic  partial   differential equations.

In the deterministic case,  in order to establish a pointwise weighted identity of the operator $L$ in \cite{F},  the operator $(\a+i\b)w_t$ was divided into $\a w_t$ and $i\b w_t$.  Then  the product of them  was estimated.  However, in the stochastic case,  the method does not work.  Therefore,  in this paper we adopt a new way to prove our pointwise weighted identity for stochastic partial differential operators, different from that in \cite{F}.

The rest of this paper is organized as follows. In Section
\ref{ss2},   a pointwise weighted identity  for some stochastic partial differential operators is established.  Section \ref{ss3} is devoted to its applications in control problems for deterministic/stochastic partial differential equations.  As  its another application, in Section \ref{ss4},   an inverse problem for linear stochastic complex Ginzburg-Landau equations  is studied. Finally,
Appendix A is  given to prove a Carleman estimate for stochastic heat equations.

\section{A  pointwise weighted identity for stochastic partial differential operators }\label{ss2}

Let $T>0$ and $\left(\Omega, \mathcal{F}, \{\mathcal{F}_t\}_{t\geq 0},
\mathcal{P}\right)$ be a complete filtered probability space, on
which a one-dimensional standard Brownian motion $\{B(t)\}_{t\geq
0}$ is defined such that $\dbF=\{\mathcal{F}_t\}_{t\geq 0}$ is the
 natural filtration generated by $B(\cdot)$, augmented by all the $\mathcal{P}$-null sets in $\mathcal{F}$. Also,  for any complex number $c$,  we denote  by $\overline c$, $\Re c$ and $\Im c$,  its complex conjugate,  real part and imaginary part, respectively.

For any  $a_0,  a,  b\in\dbR$,   $a^{j k}=a^{k j}\in L^2_{\dbF}(\Omega; C^1([0, T]; W^{2, \infty}(\dbR^n; \dbR)))$ $(j, k=1, \cdots, n)$ and ${\bf b_0}=(b^1_0, \cdots, b^n_0)\in\dbR^n$, we define the following  complex stochastic partial differential   operator:
\be\label{cp}
 \cL w=a_0dw-(a+ib)\sum\limits_{j, k=1}^{n}(a^{j k}w_{x_j})_{x_k}dt+{\bf b_0}\cdot \n wdt.
 \ee
This section is devoted to establishing a pointwise weighted identity for the operator $\cL$.  To begin with, we introduce the following assumptions:
 \begin{enumerate}  \item  If $ a\neq0$,   then ${\bf b_0}\cdot \n w$  is a lower order term for   the operator $\cL$. In this case,  without loss of generality, we  assume that ${\bf b_0}={\bf 0}$.

\item    If $a=0$ and $a_0, b\neq0$,   then  $\cL$ is a second-order stochastic Schr\"odinger operator. In this case,  we assume that ${\bf b_0}={\bf 0}$.

\item   If $a=b=0$, we assume that $a_0\neq 0$ and  ${\bf b_0}\neq{\bf 0}$. In this case,  $\cL$ is a first order stochastic transport operator.
\end{enumerate}

  The main result of this paper is stated as follows.
 \bt\label{t1}  Assume that the assumptions $1$-$3$ hold.  Let $\ell\in C^3(\dbR^{n+1}; \dbR)$, $\Phi\in C^1(\dbR^{n+1};\ \dbC)$ and $w$ be an $H^2(\dbR^n; \ \dbC)$-valued continuous semimartingle. Set $\theta=e^\ell$ and $z=\theta w$. Then for a.e. $(x, t)\in \dbR^{n+1}$ and $\mathcal{P}$-a.s. $\omega\in\Omega$,  one has the following pointwise weighted identity:
\begin{eqnarray}\label{a1}\begin{array}{rl}&2\Re(\th\overline{I_1}\mathcal{L}w)\\
\ns&\ds=2|I_1|^2dt+dM+\sum\limits_{k=1}^{n} V^k_{x_k}+B|z|^2dt+\sum\limits_{j, k=1}^n D^{j k} z_{x_j}\overline{z}_{x_k}dt\\
 \ns&\ds\q +2\sum_{j=1}^n\Big\{\Re \Big[\Big(a E^j+\overline{\Phi}  b_0^j\Big)  \overline {z}z_{x_j}\Big]+b\Im\left(F^j z\overline{z}_{x_j}\right)\Big\}dt-aa_0\sum\limits_{j, k=1}^n  a^{j k}dz_{x_j}d\overline{z}_{x_k}\\
 \ns&\ds\q-{\bf b_0}\cdot\n\[(a_0\ell_t+{\bf b_0}\cdot\n\ell)|z|^2\]dt +a_0(aA+a_0\ell_t+{\bf b_0}\cdot\n\ell)|dz|^2\\
\ns&\q
\ds-2a_0b\sum\limits_{j, k=1}^{n}a^{j k}\ell_{x_k}\Im (dz d\overline{z}_{x_j})+2 a_0\Big[b \sum\limits_{j, k=1}^{n}(a^{j k}\ell_{x_k})_{x_j}\Im (zd\overline{z})+\Re(\overline{\Phi}\overline{z}dz)\Big],
\end{array}\end{eqnarray}
where
\begin{eqnarray}\label{a2}
\left\{\begin{array}{rl}
&A=\ds\sum\limits_{j, k=1}^n \Big[a^{j k}\ell_{x_j}\ell_{x_k}-(a^{j k}\ell_{x_j})_{x_k}\Big],\quad \L=\sum\limits_{j, k=1}^n (a^{j k}z_{x_j})_{x_k}+Az,\\ \ns
&I_1=\ds-a\L+2ib\sum\limits_{j, k=1}^n a^{j k}\ell_{x_j}z_{x_k}+(\Phi-a_0\ell_t-{\bf b_0}\cdot\n\ell)z,
\end{array} \right.\end{eqnarray}
and
\begin{eqnarray*}
\left\{
\begin{array}{rl}
&B=\ds2(a^2+b^2)\sum\limits_{j, k=1}^n (A a^{j k}\ell_{x_j})_{x_k} +aa_0A_t+2aA \Re\Phi-2bA\Im \Phi\\
\ns
&\ds\quad\quad-2\Re\[\Phi(\overline{\Phi}-a_0\ell_t-{\bf b_0}\cdot\n\ell)\]+a_0\[a_0\ell_{tt}+({\bf b_0}\cdot\n\ell)_t\]+{\bf b_0}\cdot\n(a_0\ell_t+{\bf b_0}\cdot\n\ell),\\ \ns
&D^{j k}=-aa_0 a^{j k}_t +2b\Im \Phi a^{j k}-2a\Re\Phi a^{j k}\\
\ns
&\ds\quad\quad\quad+2(a^2+b^2)\sum\limits_{j', k'=1}^{n}\[a^{j k'}(a^{j'k}\ell_{x_{j'}})_{x_{k'}}+a^{kk'}(a^{j'j}\ell_{x_{j'}})_{x_{k'}}
-(a^{j k}a^{j' k'}\ell_{x_{j'}})_{x_{k'}}\],  \\ \ns
&M=\ds-aa_0A|z|^2+a_0\sum\limits_{j, k=1}^{n} a^{j k}\Big[az_{x_j}\overline{z}_{x_k}
+2b\ell_{x_j}\Im(\overline{z}_{x_k}z)\Big]-a_0(a_0\ell_t+{\bf b_0}\cdot\n\ell)|z|^2,\\
\ns
&V^k=\ds-2aa_0\sum\limits_{j=1}^{n}a^{j k}\Re(z_{x_j}d\overline{z})
-2a_0b\sum\limits_{j=1}^{n}a^{j k}\ell_{x_j}\Im(z d\overline{z})-2A(a^2+b^2)\sum\limits_{j=1}^{n} a^{j k}\ell_{x_j}|z|^2dt\\
 \ns
&\ds\quad\quad\quad
+2a\sum\limits_{j=1}^{n} a^{j k}\Re(\overline{z}_{x_j}{\Phi}z)dt+2b\sum\limits_{j=1}^{n}a^{j k}\Im\[z_{x_j}(\overline{\Phi}-a_0\ell_t)\overline{z}\]dt\\ \ns
&\ds\quad\quad\quad+2(a^2+b^2)\sum\limits_{j, j', k'=1}^{n} \[a^{jk}a^{j'k'}\ell_{x_j}z_{x_{j'}}{\overline {z}}_{x_{k'}}-a^{jk'}a^{j'k}\ell_{x_j}(z_{x_{j'}}{\overline{z}}_{x_{k'}}+{\overline{z}}_{x_{j'}}z_{x_{k'}})\]dt,\\
\ns
 &\ds E^j=\sum_{k=1}^n a^{j k}\[2\ell_{x_k}(\overline{\Phi}-a_0\ell_{t})-\overline{\Phi}_{x_k}\], \\  \ns
 &\ds F^j=\sum\limits_{k=1}^n
 \[a^{j k}(\Phi-a_0\ell_t)_{x_k}-a_0(a^{j k}\ell_{x_k})_t-2a^{j k}\ell_{x_k}\Phi\].
\end{array}
\right.
\end{eqnarray*}
\et

 \br
The pointwise weighted identity $(\ref{a1})$ is quite useful in deriving global Carleman estimates for the deterministic/stochastic partial  differential operator $(\ref{cp})$.  The advantage of Carleman inequalities derived by the identity $(\ref{a1})$  is that one can give an explicit estimate on constants  $($in Carleman estimates$)$.  This is crucial in studying  nonlinear controllability and observability  problems.\er
\br
The key point of proving the identity $(\ref{a1})$ is to multiply   ``the principal operator $\cL$"  by a weighted  multiplier $\th\overline {I_1}$.  One can  rewrite this product as a sum of ``divergence" terms,  ``energy"  terms and some lower order terms. Also, all terms in the right side of the sign of equality in $(\ref{a1})$ are real-valued functions. By choosing a suitable auxiliary function $\Phi$ and a weighted function $\th$, one can derive  global Carleman estimates for  some deterministic/stochastic partial differential  operators of different type.
 \er

\br
If choosing different coefficients in $(\ref{cp})$, one can get   deterministic/stochastic partial differential operators of different type. For example,  suppose that     $(a^{j k})_{1\leq j, k\leq n}$ is a uniformly positive definite matrix and $a_0=1$.  If $a=0$ and $b\neq 0$,  $\cL$ is a  stochastic Schr\"odinger  operator.  If $a\neq 0$,  $\cL$ is a linear stochastic complex Ginzburg-Landau operator.
If $a\neq 0$,  $b=0$ and all functions are real-valued,  $\cL$ is a stochastic heat operator.    If $a=b=0$ and ${\bf b_0}\neq {\bf0}$,
$\cL$ is a stochastic transport operator.  Also, if all functions $($in the above operators$)$ are
independent of sample points, then one can get a deterministic  Schr\"odinger operator,  a deterministic linear complex Ginzburg-Landau operator, a deterministic heat operator and a deterministic transport operator, respectively.
 In the following sections, we use the  pointwise weighted identity $(\ref{a1})$ to derive global  Carleman estimates for the above
deterministic/stochastic partial differential operators. Moreover, it is applied to study inverse problems of  linear stochastic complex Ginzburg-Landau equations.
\er

\smallskip

\noindent{\bf  Proof of Theorem \ref{t1}. }  The whole proof is divided   into four steps.

\ms

\noindent{\bf Step 1. }  Set $\theta=e^{\ell}$ and $z=\theta w$. Then it is easy to show that
\begin{eqnarray*}
&&\theta \mathcal{L}w=a_0\theta d(\theta^{-1}z)-\theta(a+ib)\sum\limits_{j, k=1}^{n}
[a^{j k}(\theta^{-1}z)_{x_j}]_{x_k}dt+\th {\bf b_0}\cdot\n (\theta^{-1}z)dt=I_1dt+I_2,
\end{eqnarray*}
where $I_1$ is given in  (\ref{a2}) and
\begin{eqnarray*}
I_2=a_0dz-ib\L dt+2a\sum\limits_{j, k=1}^{n}a^{j k}\ell_{x_j}z_{x_k}dt+{\bf b_0}\cdot\n zdt-\Phi zdt.
\end{eqnarray*}
Therefore,
\begin{equation}\label{2}
2\Re(\th\overline{I_1}\mathcal{L}w)=\theta(\overline{I_1}\mathcal{L}w+
I_1\overline{\mathcal{L}w})
=2|I_1|^2dt+2\Re(\overline{I_1}I_2).
\end{equation}

\medskip

\noindent{\bf Step 2. }  Let us compute   ``$2\Re(\overline{I_1}I_2)$".  By the assumptions 1-3, it is easy to find that
 \bel{hk1}
 a{\bf b_0}=0\quad \mbox{and} \quad  b{\bf b_0}=0.
 \ee
Recalling the definitions of $I_1$ and $I_2$,  by (\ref{hk1}) and a short calculation, we have  that

\begin{eqnarray}\label{03}
\begin{array}{ll}
2\Re(\overline{I_1}I_2)\\
\ns\ds=-2aa_0\Re(\overline{\L}dz)-4(a^2+b^2)\Re\sum_{j,k=1}^na^{jk}\ell_{x_j}(z_{x_k}\overline{\L})dt+2a\Re(\Phi\overline {\L} z)dt\\
\ns\ds\q+4a_0b\sum\limits_{j,k=1}^na^{jk}\ell_{x_j}\Im(\overline{z}_{x_k}dz)
+4b\sum\limits_{j,k=1}^na^{jk}\ell_{x_j}\Im(\overline{\Phi}\overline zz_{x_k})dt\\
\ns\ds\q+2b\Im\[(\overline{\Phi}-a_0\ell_t)\overline{z}\L\]dt+4a\sum\limits_{j,k=1}^na^{jk}\ell_{x_j}\Re\[(\overline{\Phi}-a_0\ell_t)\overline{z}z_{x_k}\]dt\\
\ns\ds\q+2\Re\[(\overline{\Phi}-a_0\ell_t-{\bf b_0}\cdot\n\ell)\overline{z}(a_0dz+{\bf b_0}\cdot\n zdt)\]\\
\ns\ds\q-2\Re\[\Phi(\overline{\Phi}-a_0\ell_t-{\bf b_0}\cdot\n\ell)\]|z|^2dt.
 \end{array}
 \end{eqnarray}

 \medskip

 \noindent{\bf Step 3. } Now we compute every term in the right side of the  sign of equality  in (\ref{03}), respectively.
By  (\ref{a2}), we find that
\begin{eqnarray}\label{3}
\begin{array}{rl} & -2aa_0\Re(\overline{\L}dz)\ds=-aa_0(\overline{\L}dz+\L d{\overline z})\\  \ns&\ds=-aa_0\sum\limits_{j, k=1}^{n}\[(a^{j k}z_{x_j})_{x_k}d\overline{z}+(a^{j k}\overline{z}_{x_j})_{x_k}dz\]-aa_0A( z d\overline{z}+\overline {z}d z)\\
\ns
&\ds
=-aa_0\sum\limits_{j, k=1}^{n}(a^{j k}z_{x_j}d\overline{z}+a^{j k}\overline{z}_{x_j}dz)_{x_k}+\sum\limits_{j, k=1}^{n} d(aa_0a^{j k}z_{x_j}\overline{z}_{x_k})-aa_0\sum\limits_{j, k=1}^{n}a^{j k}_t z_{x_j}\overline{z}_{x_k}dt\\
\ns&\ds\quad-aa_0\sum\limits_{j, k=1}^{n}a^{j k}dz_{x_j}d\overline{z}_{x_k}-d(aa_0A|z|^2)+aa_0A_t|z|^2dt+aa_0A|dz|^2.
\end{array}
 \end{eqnarray}
Further,
\begin{eqnarray}\label{5}
\begin{array}{rl}&\ds
-4(a^2+b^2)\Re\sum_{j,k=1}^na^{jk}\ell_{x_j}(z_{x_k}\overline{\L})dt=-2(a^2+b^2)\sum\limits_{j,k=1}^na^{jk}\ell_{x_j}(\overline{z}_{x_k}\L+z_{x_k}\overline{\L})dt\\
\ns
&\ds=-2(a^2+b^2)\sum\limits_{j, k=1}^{n} a^{j k}\ell_{x_j} (\overline{z}_{x_k}Az+z_{x_k}A\overline{z})dt\\ \ns
&\quad\ds-2(a^2+b^2)\sum_{j ,k=1}^{n}a^{j k}\ell_{x_j}\left[\overline{z}_{x_k}
\sum_{j', k'=1}^n(a^{j' k'}z_{x_{j'}})_{x_{k'}}+
z_{x_k}\sum_{j', k'=1}^n(a^{j' k'}\overline{z}_{x_{j'}})_{x_{k'}}\right]dt \\ \ns
&=\ds-2(a^2+b^2)\sum\limits_{j,k=1}^n(Aa^{jk}\ell_{x_j}|z|^2)_{x_k}dt+2(a^2+b^2)\sum\limits_{j,k=1}^n(Aa^{jk}\ell_{x_j})_{x_k}|z|^2dt\\ \ns
&\ds\quad-2(a^2+b^2)\sum\limits_{j, k, j', k'=1}^{n}\[
a^{j k}\ell_{x_j}a^{j' k'}(z_{x_{j'}}\overline{z}_{x_{k}}+\overline{z}_{x_{j'}}z_{x_k})\]_{x_{k'}}dt\\
\ns&\ds\q+2(a^2+b^2)\sum\limits_{j, k, j', k'=1}^{n}a^{j' k'}(a^{j k}\ell_{x_j})_{x_{k'}}(z_{x_{j'}}\overline{z}_{x_{k}}+\overline{z}_{x_{j'}}z_{x_k})dt\\
\ns&\ds\quad+2(a^2+b^2)\sum\limits_{j, k, j', k'=1}^{n} \left[(
a^{j k}\ell_{x_j}a^{j' k'}z_{x_{j'}}\overline{z}_{x_{k'}})_{x_k}-(
a^{j k}a^{j' k'}\ell_{x_j})_{x_k}z_{x_{j'}}\overline{z}_{x_{k'}}\right]dt.
  \end{array}
 \end{eqnarray}
Notice that in the above derivation,  we  use the following identity:
\begin{eqnarray*}
\begin{array}{rl}
&\ds2\sum\limits_{j, k, j', k'=1}^{n}a^{j k}a^{j' k'}\ell_{x_{j}} (z_{x_{j'}}\overline{z}_{x_kx_{k'}}+\overline{z}_{x_{j'}}z_{x_kx_{k'}})dt\\
\ns
&=\displaystyle \sum\limits_{j, k, j', k'=1}^{n}\Big\{\Big[
a^{j k}a^{j' k'}\ell_{x_{j}}(z_{x_{j'}}\overline{z}_{x_{k'}}+
\overline{z}_{x_{j'}}z_{x_{k'}})\Big]_{x_{k}} -(a^{j k}a^{j' k'}\ell_{x_{j}})_{x_{k}}(z_{x_{j'}}\overline{z}_{x_{k'}}+\overline{z}_{x_{j'}}z_{x_{k'}})\Big\}dt\\\ns
&\ds=2\sum\limits_{j, k, j', k'=1}^{n} \left[(
a^{j k}a^{j' k'}\ell_{x_j}z_{x_{j'}}\overline{z}_{x_{k'}})_{x_k}-(
a^{j k}a^{j' k'}\ell_{x_j})_{x_k}z_{x_{j'}}\overline{z}_{x_{k'}}\right]dt.
\end{array}
\end{eqnarray*}
Further,
 \begin{eqnarray}\label{6}
  \begin{array}{rl}
 &2a\Re(\Phi\overline {\L} z)dt\ds=2a\sum\limits_{j, k=1}^{n}\Re\[(a^{j k}\overline{z}_{x_j})_{x_k}\Phi z\]dt+2aA\Re \Phi |z|^2dt\\
\ns
&\ds=2a\sum\limits_{j, k=1}^{n}\Re\left(a^{j k}\overline{z}_{x_j}\Phi z\right)_{x_k}dt-2a \Re\Phi\sum\limits_{j,k=1}^n a^{j k}z_{x_j}\overline{z}_{x_k}dt\\
\ns&\ds\quad-2a\sum\limits_{j,k=1}^n \Re\left(a^{j k}\Phi_{x_k}z\overline{z}_{x_j}\right)dt+2aA \Re\Phi|z|^2dt.
\end{array}
\end{eqnarray}
Note that for  any $k=1, \cdots, n$,
 \bel{c02}
\Im(\overline{z}_{x_k}dz)=\Im\[d(\overline{z}_{x_k}z)-(zd\overline z)_{x_k}-d\overline{z}_{x_k}dz+z_{x_k}d\overline z\] =-\Im (z_{x_k}d\overline z).\ee
Therefore, we get that
\begin{eqnarray}\label{10}
\begin{array}{rl}&\ds4a_0b\sum\limits_{j,k=1}^na^{jk}\ell_{x_j}\Im(\overline{z}_{x_k}dz)=2a_0b\sum\limits_{j, k=1}^{n}a^{j k}\ell_{x_j}\Im\[d(\overline{z}_{x_k}z)-(zd\overline z)_{x_k}-d\overline{z}_{x_k}dz\]\\
\ns&\ds=2a_0b\sum\limits_{j, k=1}^{n}\Big\{d\Big[a^{j k}\ell_{x_j}\Im(\overline{z}_{x_k}z)\Big]-\Big[a^{j k}\ell_{x_j}\Im(zd\overline{z})\Big]_{x_k}\Big\}\\
\ns&\ds\quad-2a_0b\sum\limits_{j, k=1}^{n}\[(a^{j k}\ell_{x_j})_t\Im (\overline{z}_{x_k}z)dt-(a^{j k}\ell_{x_j})_{x_k}\Im(zd\overline{z})+a^{j k}\ell_{x_j}\Im(dzd\overline{z}_{x_k})\].
\end{array}
\end{eqnarray}
Further,
 \bel{02}\ba{ll}
&2b\Im\[(\overline{\Phi}-a_0\ell_t)\overline{z}\L\]dt\ds=2b\sum\limits_{j, k=1}^n\Im \[
(a^{j k}z_{x_j})_{x_k}(\overline{\Phi}-a_0\ell_t)\overline{z}\]dt -2bA\Im {\Phi}|z|^2dt\\
\ns&\ds= 2b\sum\limits_{j, k=1}^n\Im\Big[a^{j k}z_{x_j}(\overline{\Phi}-a_0\ell_t)\overline{z}\Big]_{x_k}dt
+2b\Im {\Phi}\sum\limits_{j, k=1}^n
a^{j k}z_{x_j}\overline{z}_{x_k}dt\\
\ns&\ds\quad-2b\sum\limits_{j, k=1}^n a^{j k}\Im\[(\overline{\Phi}-a_0\ell_t)_{x_k}z_{x_j}\overline{z}\]dt-2bA\Im\Phi |z|^2dt.
\ea
\ee

\medskip
\noindent{\bf Step 4. }  Let us compute ``$2\Re\[(\overline{\Phi}-a_0\ell_t-{\bf b_0}\cdot\n\ell)\overline{z}(a_0dz+{\bf b_0}\cdot\n zdt)\]$".  Notice that
 \bel{zr2}\ba{ll}
2\Re\[(\overline{\Phi}-a_0\ell_t-{\bf b_0}\cdot\n\ell)\overline{z}(a_0dz+{\bf b_0}\cdot\n zdt)\]\\
\ns\ds=2\Re\[\overline{\Phi}\overline{z}(a_0dz+{\bf b_0}\cdot\n zdt)\]\\
\ns\ds\quad-(a_0\ell_t+{\bf b_0}\cdot\n\ell)\[a_0d(|z|^2)-a_0|dz|^2+{\bf b_0}\cdot\n (|z|^2)dt\].
\ea\ee
Further,
 \bel{zr0}\ba{ll}\ds
-(a_0\ell_t+{\bf b_0}\cdot\n\ell)\[a_0d(|z|^2)-a_0|dz|^2+{\bf b_0}\cdot\n (|z|^2)dt\]\\
 \ns\ds=-d\[a_0(a_0\ell_t+{\bf b_0}\cdot\n\ell)|z|^2\]+a_0\[a_0\ell_{tt}+({\bf b_0}\cdot\n\ell)_t\]|z|^2dt\\
 \ns\ds\quad+a_0(a_0\ell_t+{\bf b_0}\cdot\n\ell)|dz|^2\\
 \ns\ds\q-{\bf b_0}\cdot\n\[(a_0\ell_t+{\bf b_0}\cdot\n\ell)|z|^2\]dt+{\bf b_0}\cdot\n\[(a_0\ell_t+{\bf b_0}\cdot\n\ell)\]|z|^2dt.
 \ea\ee

Combining  (\ref{03})-(\ref{zr0}) with (\ref{2}), we can get the desired identity (\ref{a1}).\endpf

\smallskip

\section{Applications in control problems for some deterministic/stochastic partial differential equations}\label{ss3}

 In this section,  we give  some concrete applications of  Theorem \ref{t1}  in deriving some known global Carleman estimates  for some deterministic/stochastic partial differential equations.
 Based on these  estimates,  one can study the controllability and observability of  deterministic/stochastic partial differential equations.

\subsection{A pointwise weighted identity for  deterministic partial differential operators}

In \cite{F},  a pointwise weighted identity  was established for the following deterministic  partial
differential operator:
$$\ds L=(\a+i\b)\pa_t+\sum\limits_{j,k=1}^n\pa_{x_k}(a^{jk}\pa_{x_j}),$$ with two real-valued
functions $\a$ and $\b$.  Based on this identity,    a universal approach of proving
  Carleman estimates  was established  to deduce  the controllability/observability results for  parabolic equations, hyperbolic equations,
Schr\"odinger equations, plate equations and linear  complex Ginzburg-Landau equations.

In this subsection, starting from Theorem \ref{t1},  one can obtain the known  weighted identity for   deterministic partial differential operators in \cite{F}. Indeed,
as a consequence of Theorem \ref{t1}, we have the following pointwise weighted identity.

\bc\label{0co1}
Suppose that  $a^{jk}=a^{kj}\in C^2(\dbR^n)$  $(j, k=1, 2, \cdots,n)$, $\ell\in C^3(\dbR^{n})$,  $\Phi\in C^1(\dbR^{n})$,  $y\in C^2(\dbR^n)$ and all functions in $(\ref{cp})$  are real-valued. Set $a_0=b=0, a=-1, {\bf b_0}={\bf 0}$, $\theta=e^\ell$ and $\ z=\theta y$. Then
\begin{eqnarray}\label{01ca1}\begin{array}{rl}
&2 \th I_1 \ds\sum\limits_{j, k=1}^{n}(a^{j k}y_{x_j})_{x_k}=\ds2|I_1|^2+\sum\limits_{k=1}^{n} V^k_{x_k}+B|z|^2+\sum\limits_{j, k=1}^n D^{j k} z_{x_j}z_{x_k}\\
&\quad\quad\quad\quad\quad\quad\quad\quad\quad-2\displaystyle\sum_{j,k=1}^n a^{j k}\(2\ell_{x_k} \Phi-\Phi_{x_k}\)zz_{x_j},
\end{array}\end{eqnarray}
where $I_1=\L+\Phi z$ with $\L$ being given by $(\ref{a2})$,
 $$\ba{ll}\ds
& V^k=\ds-2A\sum\limits_{j=1}^{n} a^{j k}\ell_{x_j}|z|^2
-2{\Phi}z\sum\limits_{j=1}^{n} a^{j k}\\
&\quad\quad\quad+2\displaystyle\sum\limits_{j, j', k'=1}^{n} \[a^{jk}a^{j'k'}\ell_{x_j}z_{x_{j'}}z_{x_{k'}}-a^{jk'}a^{j'k}\ell_{x_j}(z_{x_{j'}}z_{x_{k'}}+z_{x_{j'}}z_{x_{k'}})\],
 \ea$$
and
 \begin{eqnarray*}\left\{
\begin{array}{rl}&B=\ds2\sum\limits_{j, k=1}^n (A a^{j k}\ell_{x_j})_{x_k}-2A \Phi-2|\Phi|^2,\\ \ns

&D^{j k}=2\Phi a^{j k}+2\ds\sum\limits_{j', k'=1}^{n}\[2a^{j k'}(a^{j'k}\ell_{x_{j'}})_{x_{k'}}-(a^{j k}a^{j' k'}\ell_{x_{j'}})_{x_{k'}}\].
\end{array}\right.
\end{eqnarray*}
\ec

\medskip

If we  choose  $\a=\b=0$ in \cite {F},   the identity (\ref{01ca1}) in Corollary  \ref{0co1} is exactly the same as  \cite[Theorem 2.1]{F}. Meanwhile,  Corollary \ref{0co1} is very similar to \cite[Theorem 4.1]{FYZ}, both imply  weighted identities for elliptic operators.  Comparing  two results to each other,  we find that the only difference is  low-order terms. Therefore, it does not influence the derivation of  global Carleman estimates for elliptic operators of second order. Furthermore, by Corollary  \ref{0co1}, one also can deduce a weighted identity for deterministic hyperbolic  operators.

\medskip

If    $a_0\neq0$ and ${\bf b_0=0}$,  Theorem \ref{t1} implies a weighted identity for  deterministic parabolic operators,   deterministic Schr\"odingier operator and  deterministic linear  complex Ginzburg-Landau operators, respectively.

\medskip

If   $a=b=0$, $a_0\neq0$ and  ${\bf b_0\neq0}$,  Theorem $\ref{t1}$ implies a     weighted identity for deterministic transport operators.

\medskip


\subsection{A pointwise weighted identity for   stochastic transport  operators }

In Theorem \ref{t1},  assume that all functions are real-valued. If we choose $a=b=0$ and $a_0=1$,  then we have  the following  weighted identity for the stochastic transport operator: $$\mathcal{L}w=dw+{\bf b_0}\cdot\n w dt.$$

Define $\ds S^{n-1}\=\{x\in\dbR^n: |x|_{\dbR^n=1}\}$. Then, we have the following pointwise weighted identity for the stochastic trasport operator.

\bc\label{trcol}
Suppose that $\ell\in C^3(\dbR^{n+1})$ and $y$ is an $H^1(\dbR^{n})\times L^2(S^{n-1})$-valued continuous semi-martingale. Set $a=b=\Phi=0, a_0=1, {\bf b_0}\neq{\bf 0}$, $\theta=e^\ell$ and $z=\theta y$. Then
\begin{eqnarray}\label{tranport}\begin{array}{rl}
2 \th I_1 (dy+{\bf b_0}\cdot\n ydt)&\ds=2|I_1|^2dt-d\[(\ell_t+{\bf b_0}\cdot\n\ell)|z|^2\]+B|z|^2dt\\

 \ns&\ds\q -{\bf b_0}\cdot\n\[(\ell_t+{\bf b_0}\cdot\n\ell)|z|^2\]+(\ell_t+{\bf b_0}\cdot\n\ell)|dz|^2,
\end{array}\end{eqnarray}
where
\begin{eqnarray*}\left\{\ba{ll}
 I_1=\ds-(\ell_t+{\bf b_0}\cdot\n\ell)z,\\
 \ns\ds
B=\ell_{tt}+({\bf b_0}\cdot\n\ell)_t+{\bf b_0}\cdot\n(\ell_t+{\bf b_0}\cdot\n\ell).
 \end{array}\right.
\end{eqnarray*}
\ec

Corollary \ref{trcol} is exactly the same as  \cite[Proposition 2.1]{Lu3}.  As we seen in \cite{Lu3}, the identity (\ref{trcol}) plays  a key role in the study of observability/controllability problems for  stochastic transport equations.


\subsection{A global Carleman estimate for backward stochastic heat operators }

As another application of Theorem \ref{t1},  one can obtain global Carleman estimates for general forward and backward linear stochastic parabolic operators.  For simplicity, in this subsection,  we only consider backward  stochastic heat operators. Notice that our pointwise weighted identity  is different  from \cite[Theorem 3.1]{TZ}. But,   starting from this identity,  we still can obtain the desired global Carleman estimate for backward  stochastic heat equations  (which was presented in \cite{TZ}).

Let $G$ be a nonempty bounded domain in $\dbR^n$  with a boundary $\Gamma$ of class $C^4$. Put $Q=G\times(0,T)$ and $\Sigma=\Gamma\times(0,T)$.  Assume that all functions are real-valued in this subsection.  Based on the  identity  $(\ref{a1})$, we derive a global Carleman estimate for the  following    backward stochastic heat equation:
\begin{eqnarray}\label{fl000}\left\{
\begin{array}{lll}
&dy+\Delta ydt=fdt+Ydw(t)  &\mbox{ in }Q,\\[3mm]
&y=0 &\mbox{ on }\Sigma,\\[3mm]
&y(T)=y_T &\mbox{ in }G,
\end{array}
\right.
\end{eqnarray}
where $f\in L^2_\dbF(0, T; L^2(G))$ and  $y_T\in
L^2(\Omega, \mathcal{F}_T, \mathcal{P}; L^2(G))$.

\medskip

First,   introduce
some auxiliary functions.
It is well known that (\cite{fur}), there exists a   function $\psi\in C^4(\overline{G})$
such that
$$ \psi(x)>0, \mbox{ in }G;\
\ \psi(x)=0, \mbox{ on }\Gamma; \mbox{ and }  |\nabla\psi(x)|>0, \mbox{ in
}\overline{G\setminus G_1},
$$
where  $G_0$ and $G_1$ are two any given nonempty open  subsets of $G$ such that
$\overline{G_1}\subseteq G_0$. For any fixed integer $k\geq
1$, and positive parameters $\mu$ and $\lambda$, write
$$
\gamma(t)=\displaystyle\frac{1}{t^k(T-t)^k},\quad
\varphi(x,t)=\displaystyle\frac{e^{\mu\psi(x)}}{t^k(T-t)^k},\quad
\alpha(x,t)=\displaystyle \frac{e^{\mu\psi(x)}-
e^{2\mu|\psi|_{C(\overline{G})}}}{t^k(T-t)^k}\quad\mbox{and}\quad
\theta=e^{\lambda\alpha}.
$$
In the sequel,  for any $k\in\dbN$, we denote by $\mathcal{O}(\mu^k)$ a function of
 order $\mu^k$, for sufficiently large $\mu$; and by $\mathcal{O}_\mu(\lambda^k)$ a function
 of order $\lambda^k$ for fixed $\mu$ and sufficiently large $\lambda$.

\medskip

Next,  based on (\ref{a1}), we  have the following  inequality for the backward stochastic heat operator:
$$
\mathcal{L}y=dy+\Delta ydt.
$$

\begin{lemma}\label{fll1}
Let $z=\theta y$  and $\ell=\lambda\alpha$. Then any solution $(y,
Y)\in L^2_{\dbF}(0,T;H^1_0(G))\times
L^2_{\dbF}(0,T;L^2(G))$ of  the equation $(\ref{fl000})$ satisfies
\begin{eqnarray}\label{*}
\begin{array}{rl}
&\dbE\displaystyle\int_Q 2\theta I_1 \mathcal{L}y dx\\[3mm]
&\geq \dbE\displaystyle\int_Q 2|I_1|^2dxdt\\[3mm]
&\quad+\dbE\displaystyle\int_Q 2\lambda^3\mu^4\varphi^3 |\nabla\psi|^4 z^2dxdt+\dbE\displaystyle\int_Q 2\lambda\mu^2\varphi |\nabla z|^2 |\nabla \psi|^2dxdt\\[3mm]
&\quad+\dbE\displaystyle\int_Q \left[
\lambda^3\varphi^3\mathcal{O}(\mu^3)
+\mathcal{O}_\mu(\lambda^2)\varphi^3 \right]|z|^2dxdt+\dbE\displaystyle\int_Q \mathcal{O}(\mu)\lambda\varphi|\nabla z|^2dxdt\\[3mm]
&\quad+\dbE\displaystyle\int_Q \theta^2\mathcal{O}(\lambda^2)\mu^2\varphi^2|Y|^2
dxdt,\end{array}
\end{eqnarray}
\end{lemma}
where $I_1=\Delta z+(|\nabla \ell|^2+\Delta \ell-\ell_t)z$.

\medskip

For the readers' convenience, we give a proof of Lemma \ref{fll1} in Appendix A.  By (\ref{*}),   proceeding exactly the same analysis as
 \cite[Theorem 6.1]{TZ},  one can obtain the following global Carleman estimate for the equation (\ref{fl000}).

\begin{theorem}\label{tfl11}
There exists a positive constant
$\mu_0$, depending only on $n, G,
 G_0$ and $T$,  such that for any $\mu\geq \mu_0$, one can find two positive constants
$\lambda_0=\lambda_0(\mu)$ and $C=C(\mu)$, so that for
any $f\in L^2_{\dbF}(0,T;L^2(G))$ and $y_T\in
L^2(\Omega,\mathcal{F}_T,\mathcal{P};L^2(G))$,  any solution $(y,
Y)\in L^2_{\dbF}(0,T;H^1_0(G))\times
L^2_{\dbF}(0,T;L^2(G))$ of  the equation $(\ref{fl000})$ satisfies
\begin{eqnarray*}
&&\dbE\int_Q \theta^2\left(\lambda^3 \gamma^3 y^2+
\lambda\gamma|\nabla y|^2\right) dxdt\\[3mm]
&&\leq
C\left(\dbE\int^T_0\int_{G_0}\theta^2\lambda^3\gamma^3 y^2dxdt+
\dbE\int_Q
\theta^2f^2dxdt+\dbE\int_Q\theta^2\lambda^2\gamma^2
Y^2dxdt\right),
\end{eqnarray*}
for any $\lambda\geq \lambda_0.$
\end{theorem}

By this global Carleman estimate in Theorem \ref{tfl11},  one  can study the observability (resp. controllability) for backward (resp. forward) stochastic heat equations.

\subsection{A global Carleman estimate  for stochastic Schr\"odinger  operators }

In this subsection, based on the weighted identity  (\ref{a1}),  we  derive a global Carleman estimate  for  stochastic
Schr\"odinger operators.
In  (\ref{cp}),  we choose
$
a_0=1, b=1, a=0, {\bf b_0}={\bf 0}
$ and $(a^{j k})_{1\leq j, k\leq n}=I_n$. Then $\mathcal{L}$ is the following stochastic
Schr\"odinger operator:
\begin{equation}\label{fl20}
\mathcal{L}w=dw-i\Delta wdt,
\end{equation}
and $I_1=2i\nabla \ell\cdot\nabla z+(\Phi-\ell_t)z$ (with $z=\theta w$).

\medskip

Notice that  in   \cite{Lu1}, a weighted identity was derived,  in order to establish a global Carleman estimate for the stochastic
Schr\"odinger operator:
$$
Pv=idv+\Delta vdt.
$$
In \cite{Lu1},   write $u=\theta v$. Then
$$
\widetilde{I_1}=-i\ell_t u-2\nabla\ell\cdot\nabla u+\Psi u,
$$
where $\Psi$ is a suitable auxiliary function.

If we set $w=iv$ and $\Phi=-i\Psi$,  then it is easy to check that $z=iu$, $I_1=\widetilde{I_1}$ and
$\mathcal{L}w=Pv$. Therefore,  based on (\ref{a1}), we can get  the same pointwise weighted identity as that in \cite{Lu1}, and a  global Carleman estimate for stochastic
Schr\"odinger operators.

\section{Applications in inverse problems for linear stochastic complex Ginzburg-Landau equations}\label{ss4}

As another application of Theorem \ref{t1}, in this section, we prove a  uniqueness result for inverse problems of  linear  stochastic complex Ginzburg-Landau equations.

\subsection{Main results}

The deterministic complex Ginzburg-Landau equation was introduced by Ginzburg and Landau in 1950  (\cite{GL}). This kind of complex partial differential equations can describe a phase transition in the theory of superconductivity. In the last decades, a lot of stochastic versions of Ginzburg-Landau equations were  studied. We refer to \cite{BMP, LLL} and the references therein for some known results.

Consider the following linear stochastic complex Ginzburg-Landau equation:

  \bel{0c1}\left\{\ba{ll}\ds dw-(1+ib)\D wdt=(a_1\cdot\n w+a_2w)dt+a_3wdB(t) &\mbox{ in } Q,\\
 \ns \ds w=0 &\mbox{ on }\Si,\\
  \ns \ds w(0)=w_0&\mbox{ in }G,
   \ea\right.\ee
 where $b\in\dbR$,\  $a_1\in L^\i_{\dbF}(0,T; L^\i(G;\ \dbC^n)),\ a_2\in L^\i_{\dbF}(0,T; L^\i(G;\ \dbC)),\  a_3\in L^\i_{\dbF}(0,T; W^{1,\i}(G;\ \dbC))$ and $w_0\in L^2(G; \ \dbC)$.
 \medskip

We first  recall the definition of  weak solutions of the equation (\ref{0c1}).

 \bde\label{def1}
We call $w\in L^2_{\dbF}(\O; C([0,T];L^2(G; \ \dbC)))\bigcap L^2_{\dbF}(0,T; H_0^1(G; \ \dbC))$ is a weak solution of the equation $(\ref{0c1})$,  if for any $t\in[0,T]$ and any $p\in H_0^1(G)$, it holds that
 $$\ba{ll}\ds
\int_G w(t,x)\overline{p}(x)dx-\int_G w_0\overline{p}(x)dx\\
\ns\ds=\int_0^t\int_G\Big\{-(1+ib)\nabla w(s,x)\cdot\nabla\overline{p}(x)+\[a_1(s,x)\cdot\n w(s,x)+a_2(s,x)w(s,x)\]\overline{p}(x)\Big\}dxds\\
\ns\ds\q+\int_0^t\int_G a_3(s,x)w(s,x) \overline{p}(x)dxdB(s), \qq \mathcal{P}\mbox{-}a.s. \ea $$\ede
Also, set
  \bel{rr}
  r\= 1+|a_1|^2_{L^\i_{\dbF}(0,T; L^\i(G;\ \dbC^n))}+|a_2|^2_{L^\i_{\dbF}(0,T; L^\i(G;\ \dbC))} +|a_3|^2_{L^\i_{\dbF}(0,T; W^{1,\i}(G;\ \dbC))}.
  \ee
Then we have the following well-posedness result for the equation (\ref{0c1}), whose proof can be found in \cite[Chapter 6]{PZ}.

 \bl\label{lem3}  For any $w_0\in L^2(G; \ \dbC)$, there exists  a unique weak solution $w$ of  the equation $(\ref{0c1})$ . Moreover,   $$
 |w|_{ L^2_{\dbF}(\O; C([0,T];L^2(G;\  \dbC)))}+|w|_{L^2_{\dbF}(0,T;H_0^1(G;\ \dbC))}\le C r |w_0|_{L^2(G;\  \dbC)}.
 $$   \el

In this section, we are concerned with the following inverse problem: for $t_0\in [0,T)$, determine $w(t_0,\cdot)$, $\mathcal{P}\mbox{-}a.s.$ from $w(T,\cdot)$.  The main result can be stated as follows.

 \bt\label{t4}  Let $t_0\in[0,T)$. Then there exist  constants $\tau\in(0,1)$ and  $C>0$ such that
  \bel{cc07}
  |w(t_0)|_{L^2(\O,\cF_{t_0}, \mathcal{P}; L^2(G; \ \dbC))}\le C|w|^{1-\tau}_{L^2_{\dbF}(0,T;L^2(G; \ \dbC))}|w(T)|^\tau_{L^2(\O,\cF_T, \mathcal{P}; H^1(G; \ \dbC))},
  \ee
 for any solution $w$  of the  equation $(\ref{0c1})$.
       \et

  As a consequence of Theorem \ref{t4}, we can get the following  backward uniqueness for the equation (\ref{0c1}).
   \bc\label{0c2} Assume that $w$ is a weak solution of the equation  $(\ref{0c1})$.
     If $w(T)=0$ in $G$, $\mathcal{P}$-a.s., then $w(t)=0$ in $G$, $\mathcal{P}$-a.s. for all $t\in[0,T]$.
     \ec

 In  \cite{F, RZ},
some global Carleman estimates for  deterministic complex Ginzburg-Landau equations were established, respectively.
However, as far as we know,  there are no published papers addressing  global Carleman estimates for stochastic complex Ginzburg-Landau equations.  In the following, we derive a suitable Carleman estimate for a linear stochastic complex Ginzburg-Landau oeprator. Based on this result,  we can study the  uniqueness  of this inverse problem.

\subsection{A  Carleman estimate for  linear stochastic complex Ginzburg-Landau operators }

In this subsection, we establish a  Carleman estimate for the following linear stochastic complex  Ginzburg-Landau equation:

 \bel{01c1}\left\{\ba{ll}\ds dw-(1+ib)\D wdt=fdt+gdB(t) &\mbox{ in } Q,\\
 \ns \ds w=0 &\mbox{ on }\Si,\\
  \ns \ds w(0)=w_0&\mbox{ in }G,
   \ea\right.\ee
 where $f\in L^2_{\dbF}(0,T; L^2(G;\ \dbC))$ and  $g\in L^2_{\dbF}(0,T; H^{1}(G;\ \dbC))$.

 \medskip

 First, we establish a pointwise weighted identity, which is a consequence  of Theorem \ref{t1}.
  \bl\label{lem2} Under the assumptions of Theorem $\ref{t1}$, for a parameter $\mu\geq 1$,  choose
  $\ds
\varphi(t)=e^{3\mu t},\ \ell=\mu \varphi, \th=e^{\ell}, z=\th w$ and $\Phi=-\mu.
  $
Then,  it holds that
 \begin{eqnarray}\label{2a1}\begin{array}{rl}
&2\Re\[\th\overline{I_1}\(dw-(1+ib)\D wdt\)\]\\
\ns&\ds=\ds2|I_1|^2dt+ d(|\n z|^2-3\mu^2 \varphi|z|^2)+\sum\limits_{k=1}^{n} V^k_{x_k}+\mu^2(3\mu\varphi-2)|z|^2dt\\
 &\ds\q+2\mu |\n z|^2dt-|\n dz|^2-2\mu\Re(\overline{z}dz)+3\mu^2\varphi |dz|^2,
\end{array}\end{eqnarray}
where
 \bel{cc02}\left\{\ba{ll}\ds
I_1=\ds-\D z-(\mu+3\mu^2 \varphi) z,\\
\ns\ds V^k=\ds-2\Re(z_{x_k}d\overline{z}
+\mu \overline{z}_{x_k}zdt)-2b(\mu+3\mu^2 \varphi)
\Im(z_{x_k}\overline{z})dt.
\ea\right.
 \ee
 \el

 \noindent {\bf Proof. }   In Theorem \ref{t1},   we choose $a_0=a=1, {\bf b_0=0},  (a^{jk})_{n\t n}=I_n$, $\ds
\varphi(t)=e^{3\mu t},\ \ell=\mu \varphi$ and $\Phi=-\mu$. Then after a simple calculation, we can get the desired result (\ref{2a1}). \endpf

 \medskip

\medskip

Based on Lemma \ref{lem2}, we have  the following Carleman estimate for (\ref{01c1}).

 \bt\label{t3} Let $\d\in[0,T)$. Then  for any  $\mu\ge2$, one can find a
constant $C=C(\mu)>0$ so that
 \bel{pp}\ba{ll}&\ds
\mu\dbE\int_{\d}^T\int_G\th^2|\n w|^2dxdt+\mu^3\dbE\int_{\d}^T\int_G\varphi \th^2|w|^2dxdt\\
 \ns&\ds\le
 C\Big\{\dbE \int_G\[|\th(\d)\n w(\d)|^2+\mu^2 \varphi(\d)\th(\d)|w(\d)|^2+\mu^2 \varphi(T)|\th(T)w(T)|^2\]dx\\
 \ns&\ds\qq+\dbE\int_{\d}^T\int_G(1+\varphi)\th^2(|f|^2+\mu^2|g|^2+|\n g|^2)dxdt\Big\}. \ea\ee
for any  solution $w$  of the equation $(\ref{01c1})
$.
   \et

\noindent {\bf Proof. }   Integrating the identity (\ref{2a1})  in $[\d,T]\t G$ for  $\d\in[0,T)$, and taking mathematical expectation,  by (\ref{01c1}) and  $z|_{\Si}=0$,  we have that
  \begin{eqnarray}\label{c2a1}\begin{array}{rl}
&\ds\ds2\dbE\int_{\d}^T\int_G|I_1|^2dxdt+ \dbE\int_{\d}^T\int_Gd(|\n z|^2-3\mu^2 \varphi|z|^2)dx+\mu^2\dbE\int_{\d}^T\int_G(3\mu\varphi-2)|z|^2dxdt\\
\ns&\ds\qq+2\dbE\int_{\d}^T\int_G\mu |\n z|^2dxdt-\dbE\int_{\d}^T\int_G\[|\n dz|^2+2\mu\Re(\overline{z}dz)-3\mu^2\varphi dzd\overline{z}\]dx\\
 \ns&\ds=2\dbE\int_{\d}^T\int_G\Re\[\th\overline{I_1}(fdt+gdB)\]dx\\
 \ns&\ds\le 2\dbE\int_{\d}^T\int_G|I_1|^2dtdx+2\dbE\int_{\d}^T\int_G|\th f|^2dxdt.
\end{array}\end{eqnarray}
It is easy to check that
  \bel{cc04}\ba{ll}\ds
 -\dbE\int_{\d}^T\int_Gd(|\n z|^2-3\mu^2 \varphi|z|^2)dx\ds\le  C\dbE \int_G\[|\n z(\d)|^2+\mu^2  \varphi(T)|z(T)|^2\]dx.
   \ea\ee

 On the other hand, noting that  $ 2\Re(\overline{z}dz) =zd{\overline z}+{\overline z}dz=d(|z|^2) -|dz|^2$, we obtain that

   \bel{cc05}\ba{ll}\ds
 \dbE\int_{\d}^T\int_G\[|\n dz|^2-3\mu^2\varphi dzd\overline{z}+2\mu\Re(\overline{z}dz)\]dx\\
 \ns\ds\le C\dbE\int_{\d}^T\int_G\th^2\[|\n g|^2+\mu^2(1+\varphi) \th^2|g|^2 \]dxdt+\mu\dbE \int_G|z(\d)|^2dx.
  \ea\ee
By  (\ref{c2a1})-(\ref{cc05}),  we have that
\begin{eqnarray}\label{ca01}\begin{array}{rl}
&\ds\ds\dbE\int_{\d}^T\int_G 3\mu^2\((\mu\varphi-1)|z|^2+2\mu |\n z|^2\)dxdt\\
 \ns&\ds\le C\dbE \int_G\[|\n z(\d)|^2+\mu|z(\d)|^2+\mu^2 \varphi(T)|z(T)|^2\]dx\\
 \ns&\ds\q+C\dbE\int_{\d}^T\int_G(1+\varphi)\th^2(|f|^2+\mu^2|g|^2+|\n g|^2)dxdt.

\end{array}\end{eqnarray}

Taking $\mu_0=2$ and noting that $\varphi=e^{3\mu t}>1$,  we obtain  that $\mu_0\varphi-1>\ds\frac{\mu_0}{2}$. Therefore,
   by (\ref{ca01}) and   $z=\th w$,  one can get the desired inequality (\ref{pp}). \endpf


\subsection{Proof of Theorem \ref{t4}}

This subsection is devoted to a proof of Theorem \ref{t4}. We borrow some ideas from \cite{Lu2}.

\medskip

\noindent {\bf Proof of Theorem \ref{t4}. } The proof is  divided  into two steps.

\medskip

\noindent{\bf Step 1. } For any $t_0\in (0,T)$, we choose $t_1$ and $t_2$ satisfying that $0<t_1<t_2<t_0$.   Let $\rho\in C^\i(\dbR;[0,1])$ be a function such that
 \bel{cc10}\rho=\left\{\ba{ll}\ds 1,&t\ge t_2,\\
 \ns  0, & t\le t_1.
   \ea\right.\ee
Put $h=\rho w$. Then by the equation (\ref{0c1}), $h$ satisfies  that
    \bel{0cv1}\left\{\ba{ll}\ds dh-(1+ib)\D hdt=[(a_1,\n h)+a_2h+\rho_t w]dt+a_3hdB(t) &\mbox{ in } Q,\\
 \ns \ds h=0 &\mbox{ on }\Si,\\
  \ns \ds h(0)=0&\mbox{ in }G.
   \ea\right.\ee
Applying Theorem \ref{t3} (with $\d=0$) to the equation (\ref{0cv1}),     we can find a
$\mu_1>2$ such that for any $\mu\ge\mu_1$, \bel{0cv2}\ba{ll}&\ds
\mu\dbE\int_{0}^T\int_G\th^2|\n h|^2dxdt+\mu^3\dbE\int_{0}^T\int_G\varphi \th^2|h|^2dxdt\\
 \ns&\ds\le
 C\dbE\Big\{\th^2(T)\int_G\[|\n h(T)|^2+\mu^2 \varphi(T)|h(T)|^2\]dx+\int_Q\th^2|\rho_t(t)w|^2dxdt\Big\}.
   \ea\ee
Noting that $\th=e^{\mu e^{3\mu t}}$ is an increasing function of $t$,  by  (\ref{cc10}), we have that
 \bel{0cv3}
\dbE \int_Q\th^2|\rho_t(t) w|^2dxdt\le C\dbE\int_{t_1}^{t_2}\int_G \th^2 |w|^2dxdt\le C\th^2(t_2)|h|^2_{L^2_\dbF(0,T; L^2(G;\ \dbC))}.
  \ee
 Therefore, combining (\ref{0cv2}) and (\ref{0cv3}), we get that
    \bel{0cv4}\ba{ll}&\ds
    \mu\th^2(t_0)\dbE\int_{t_0}^T\int_G |\n h|^2dxdt+\mu^3\th^2(t_0)\dbE\int_{t_0}^T\int_G\varphi |h|^2dxdt\\
 \ns&\ds\le
\mu\dbE\int_{0}^T\int_G\th^2|\n h|^2dxdt+\mu^3\dbE\int_{0}^T\int_G\varphi \th^2|h|^2dxdt\\
 \ns&\ds\le
 C\th^2(T)\dbE\int_G\[|\n h(T)|^2+\mu^2 \varphi(T)|h(T)|^2\]dx+C\th^2(t_2)|w|^2_{L^2_\dbF(0,T; L^2(G;\ \dbC))}.
   \ea\ee
 By (\ref{0cv4}) and noting that $h=\rho w$,   we obtain that
           \bel{0cv5}\ba{ll}&\ds
    \mu\dbE\int_{t_0}^T\int_G |\n h|^2dxdt+\mu^3\dbE\int_{t_0}^T\int_G\varphi |h|^2dxdt\\
 \ns&\ds\le C\th^{-2}(t_0)\th^2(t_2)|w|^2_{L^2_\dbF(0,T; L^2(G;\ \dbC))} +C\th^2(T)\dbE\int_G\[|\n w(T)|^2+\mu^2 \varphi(T)|w(T)|^2\]dx.
   \ea\ee

\noindent{\bf Step 2. }  Let us estimate ``$\ds \dbE\int_G|w(t_0)|^2dx$".

\medskip

By (\ref{0c1}) and (\ref{rr}), it is easy to check that
      \bel{0cv6}\ba{ll}&\ds
  \dbE\int_G|w(t_0)|^2dx -\dbE\int_G|w(T)|^2dx\\
  \ns&\ds= -\dbE\int_{t_0}^T\int_G [wd{\overline w}+{\overline w}dw+|dw|^2]dx\\
  \ns&\ds=2\int_{t_0}^T\int_G|\n w|^2-\dbE\int_{t_0}^T\int_G [w(a_1,\n \overline w)+\overline w(a_1,\n  w)+2a_2|w|^2+|a_3w|^2]dxdt\\
  \ns&\ds\le  C\int_{t_0}^T\int_G|\n w|^2dxdt+Cr \dbE\int_{t_0}^T\int_G|w|^2dxdt.
    \ea\ee
   Combining (\ref{0cv5})-(\ref{0cv6}),   we find that
    \bel{0cv7}\ba{ll}\ds
  \dbE\int_G|w(t_0)|^2dx&\ds\le C\th^{-2}(t_0)\th^2(t_2)|w|^2_{L^2_\dbF(0,T; L^2(G;\ \dbC))}\\
  \ns&\quad\ds+C\mu^2 \varphi(T)\th^2(T)\dbE|w(T)|^2_{H^1(G;\ \dbC)}\\
  \ns&\ds\le Ce^{-2\mu(e^{3\mu_1 t_0}-e^{3\mu_1 t_2})}|w|^2_{L^2_\dbF(0,T; L^2(G;\ \dbC))}+Ce^{2\mu e^{C\mu T}}\dbE|w(T)|^2_{H^1(G;\ \dbC)} .
    \ea\ee
Note that  $t_2<t_0$.  We choose a $\mu>1$ as a minimizer of the right hand side in  the inequality (\ref{0cv7}). Then it follows    that
$$
  \dbE\int_G|w(t_0)|^2dx\le C|w|^{1-\tau}_{L^2_{\dbF}(0,T;L^2(G;\ \dbC))}|w(T)|^\tau_{L^2(\O,\cF_T, \mathcal{P}; H^1(G; \ \dbC))},$$
with
 $$
 \tau= \frac{2(e^{3\mu_1t_0}-e^{3\mu_1t_1})}{C+2(e^{3\mu_1t_0}-e^{3\mu_1t_1})}.
 $$
 This completes the proof of Theorem \ref{t4}.\endpf

\section{Appendix A }\label{ss5}

\noindent{\bf  Proof of Lemma \ref{fll1}. }  In (\ref{cp}),  choose $a_0=1$, $a=-1,$  $b=0$, ${\bf b_0=0}$ and $(a^{j k})_{1\leq j, k\leq n}=I_n$.  Then by $(\ref{a1})$, we obtain that
\begin{eqnarray}\label{fl00}\begin{array}{rl}&2\th I_1\mathcal{L}y\\
\ns&\ds=2|I_1|^2dt+dM+\sum\limits_{k=1}^{n} V^k_{x_k}+B|z|^2dt+\sum\limits_{j, k=1}^n D^{j k} z_{x_j}z_{x_k}dt\\
 \ns&\ds\q +2\sum_{j=1}^n \Big(E^j  zz_{x_j}\Big)dt+\sum\limits_{j=1}^n |dz_{x_j}|^2+(-A+\ell_t)|dz|^2+2\Phi zdz,
\end{array}\end{eqnarray}
where
 \begin{eqnarray*}\label{enery}\left\{
\begin{array}{rl}
&A=|\nabla \ell|^2-\Delta \ell,\quad \L=\Delta z+Az,\quad I_1=\ds\L+(\Phi-\ell_t)z,\\ \ns
&B=\ds2\sum\limits_{j=1}^n (A \ell_{x_j})_{x_j} -A_t-2A \Phi-2(\Phi^2-\ell_t\Phi)+\ell_{tt},\\ \ns
&D^{j k}=2\Phi \delta^{j}_k+4\ell_{x_{j}x_{k}}-2\Delta\ell\delta^j_k\quad\mbox{with }
\delta^j_k=\left\{\begin{array}{rl}
1\quad&j=k,\\
0\quad&j\neq k,
\end{array}\right.\\  \ns
&M=\ds A|z|^2-|\nabla z|^2-\ell_t|z|^2,\\  \ns
&V^k=\ds 2z_{x_k}dz
-2A\ell_{x_k}|z|^2dt
-2z_{x_k}\Phi zdt+2|\nabla z|^2\ell_{x_k}dt-4\nabla\ell\cdot\nabla z z_{x_k}dt,\\   \ns
 &\ds E^j=2\ell_{x_j}(\Phi-\ell_{t})-\Phi_{x_j}.
\end{array}\right.
\end{eqnarray*}
Also,  set
$
\Phi=2\Delta \ell.
$
Then it is easy to check that for any $j, k=1, \cdots, n$,
\begin{eqnarray*}
&&\ell_{x_j}=\lambda\mu\varphi\psi_{x_j},\quad \ell_{x_j x_k}=\lambda\mu\varphi\psi_{x_j x_k}+\lambda\mu^2\varphi\psi_{x_j}\psi_{x_k},\quad \ell_{tt}=\mathcal{O}_{\mu}(\lambda)\varphi^3,\\[2mm]
&&A=\lambda^2\mu^2\varphi^2|\nabla \psi|^2-\lambda\mu^2\varphi |\nabla \psi|^2-\lambda\mu\varphi \Delta \psi=\lambda^2\mu^2\varphi^2|\nabla \psi|^2+\mathcal{O}(\lambda)\mu^2\varphi,\\[2mm]
&&A_t=\mathcal{O}(\lambda^2)\mu^2\varphi^3,
\end{eqnarray*}
and for any $k=1, 2, \cdots, n$,
\begin{eqnarray*}
A_{x_{k}}=2\lambda^2\mu^3\varphi^2 |\nabla \psi|^2\psi_{x_{k}}
+\mathcal{O}(\lambda^2)\mu^2\varphi^2+\mathcal{O}(\lambda)\mu^3\varphi.
\end{eqnarray*}

\medskip

In the following, we estimate every term  in the right side of  sign of equality in (\ref{fl00}).

\medskip
\noindent{\bf Step 1. } First,  notice that
\begin{eqnarray}\label{fl1}
\begin{array}{rl}
&B=2\nabla A\cdot \nabla\ell-2A\Delta \ell
-A_t+\ell_{tt}
-8(\Delta \ell)^2
+4\Delta\ell\ell_t.
\end{array}
\end{eqnarray}
Further,
\begin{eqnarray}\label{fl2}
\begin{array}{rl}
&2\nabla A\cdot \nabla\ell=2\lambda\mu\varphi\nabla A\cdot \nabla\psi\\[3mm]
&=2\displaystyle\sum\limits_{j=1}^{n}
\left[2\lambda^2\mu^3\varphi^2 |\nabla \psi|^2 \psi_{x_{j}}
+\mathcal{O}(\mu^2)\lambda^2\varphi^2+\mathcal{O}(\lambda)\mu^3\varphi\right]\cdot\lambda\mu\varphi\psi_{x_{j}}\\[3mm]
&=4
\lambda^3\mu^4\varphi^3|\nabla\psi|^4+\mathcal{O}(\mu^3)\lambda^3\varphi^3
+\mathcal{O}(\lambda^2)\mu^4\varphi^2.
\end{array}
\end{eqnarray}
Further,
\begin{eqnarray}\label{fl3}
\begin{array}{rl}
&-2A\Delta\ell=-2\Big[\lambda^2\mu^2\varphi^2 |\nabla\psi|^2+
\mathcal{O}(\lambda)\mu^2\varphi\Big]\Big[\lambda\mu^2\varphi |\nabla\psi|^2+\mathcal{O}(\lambda)\mu\varphi\Big]\\[3mm]
&=-2\lambda^3\mu^4\varphi^3|\nabla\psi|^4+\mathcal{O}(\mu^3)\lambda^3\varphi^3+\mathcal{O}(\lambda^2)\mu^4\varphi^2.
\end{array}
\end{eqnarray}
Further,
\begin{eqnarray}\label{fl4}
\begin{array}{rl}
&-8(\Delta\ell)^2+4\Delta\ell\ell_t\\[3mm]
&=\mathcal{O}(\lambda^2)\mu^4 \varphi^2+\mathcal{O}(\lambda)\mu^2\varphi\cdot \lambda e^{2\mu|\psi|_{C(\overline{G})}}\varphi^2=\mathcal{O}(\lambda^2)\mu^4\varphi^2+\mathcal{O}_{\mu}(\lambda^2)\varphi^3.
\end{array}
\end{eqnarray}

Combining  (\ref{fl2})-(\ref{fl4})  with  (\ref{fl1}),  we  get that
\begin{eqnarray}\label{fl5}
\begin{array}{rl}
&B=2 \lambda^3\mu^4\varphi^3 |\nabla\psi|^4
+\mathcal{O}(\mu^3)\lambda^3\varphi^3+
\mathcal{O}_\mu(\lambda^2)\varphi^3.
\end{array}
\end{eqnarray}

\medskip

\noindent {\bf Step 2. }  Noticing that $z=0$ on $\Sigma$,  we have that  for any $k=1, 2, \cdots, n$,
\begin{eqnarray*}
&&V^k\Big|_{\Sigma}=2\Big(|\nabla z|^2\ell_{x_k}-
2\nabla z\cdot\nabla \ell z_{x_k}\Big)dt\Big|_\Sigma\\
&&=2\left(
\lambda\mu\varphi\Big|\frac{\partial z}{\partial \nu}\Big|^2 \psi_{x_k}
-2\lambda\mu\varphi\Big|\frac{\partial z}{\partial \nu}\Big|^2\frac{\partial \psi}{\partial\nu} \nu_{k}
\right)dt\Big|_{\Sigma},
\end{eqnarray*}
where $\nu=(\nu_1, \cdots, \nu_n)$ denotes the unit outer normal vector on $\Gamma$.  Therefore,
\begin{eqnarray}\label{fl6}
\begin{array}{rl}
&\displaystyle\sum\limits_{k=1}^{n} V^k \cdot \nu_k=-
2 \lambda\mu\varphi\Big|\frac{\partial z}{\partial \nu}\Big|^2\frac{\partial \psi}{\partial \nu} \Big|_\Sigma\geq 0.
\end{array}
\end{eqnarray}

\medskip

\noindent {\bf Step 3. } By the definitions of $D^{j k}$ $(j, k=1, \cdots, n)$ and $\Phi$,  we have
\begin{eqnarray*}
&&D^{j k}=2\Delta \ell\delta^j_k+4\ell_{x_jx_k}\\
&&=2\Big(\lambda\mu\varphi\Delta\psi+\lambda\mu^2\varphi|\nabla \psi|^2\Big)\delta^j_k
+4\lambda\mu\varphi\psi_{x_j x_k}+4\lambda\mu^2\varphi\psi_{x_j}\psi_{x_k} \\
&&=\mathcal{O}(\mu)\lambda\varphi+2\lambda\mu^2\varphi|\nabla \psi|^2\delta^j_k+4\lambda\mu^2\varphi\psi_{x_j}
\psi_{x_k}.
\end{eqnarray*}
It follows that
\begin{eqnarray}\label{fl7}
\begin{array}{rl}
&\displaystyle\sum\limits_{j, k=1}^{n} D^{j k}z_{x_j}z_{x_k}=2\lambda\mu^2\varphi |\nabla \psi|^2|\nabla z|^2+\mathcal{O}(\mu)\lambda\varphi|\nabla z|^2+4\lambda\mu^2\varphi\Big|\nabla \psi\cdot\nabla z\Big|^2.
\end{array}
\end{eqnarray}

\medskip

\noindent {\bf Step 4. }  By the definitions of $E^j$  $(j=1, 2, \cdots, n)$,  we have  that
\begin{eqnarray*}
&&2\sum\limits_{j=1}^{n} \Big(E^j  zz_{x_j}\Big)dt=2\sum\limits_{j=1}^{n}\Big(4\ell_{x_j}\Delta\ell-2\ell_{x_j}\ell_t-2\Delta\ell_{x_j}\Big)zz_{x_j}dt\\
&&=4\sum\limits_{j=1}^{n}\ell_{x_j}\Delta\ell(z^2)_{x_j}dt
-2\sum\limits_{j=1}^{n}\ell_{x_j}\ell_t(z^2)_{x_j}dt
+\mathcal{O}(\lambda)\mu^3\varphi|z||\nabla z|dt\\
&&=\sum\limits_{j=1}^{n} \Big(4\ell_{x_j}\Delta\ell z^2-2\ell_{x_j}\ell_t z^2\Big)_{x_j}dt\\
&&\quad-4\sum\limits_{j=1}^{n}(\ell_{x_j}\Delta\ell)_{x_j}z^2dt
+2\sum\limits_{j=1}^{n} (\ell_{x_j}\ell_t)_{x_j}z^2dt+\mathcal{O}(\lambda)\mu^3\varphi|z||\nabla z|dt\\
&&=\sum\limits_{j=1}^{n} \Big(4\ell_{x_j}\Delta\ell z^2-2\ell_{x_j}\ell_t z^2\Big)_{x_j}dt
+\mathcal{O}(\lambda^2)\mu^4\varphi^2z^2dt
+\mathcal{O}_\mu(\lambda^2)\varphi^3z^2dt+\mathcal{O}(\lambda)\mu^3\varphi|z||\nabla z|dt.
\end{eqnarray*}
Therefore,
\begin{eqnarray}\label{fl9}
\begin{array}{rl}
&2\dbE\ds\int_Q\sum\limits_{j=1}^{n} \Big(E^j zz_{x_j}\Big)dxdt\\[3mm]
&=\dbE\ds\int_Q \[\mathcal{O}(\lambda^2)\mu^4\varphi^2+\mathcal{O}_{\mu}(\lambda^2)\varphi^3\]
|z|^2dxdt+\dbE\ds\int_Q \mathcal{O}(\lambda)\mu^3\varphi^2|z||\nabla z|dxdt.
\end{array}
\end{eqnarray}

\medskip

\noindent{\bf Step 5. }  By the first equation of (\ref{fl000}),  we find that
\begin{eqnarray*}
&&\dbE\int_Q  (-A+\ell_t)|dz|^2 dx
=\dbE\int_Q\[\mathcal{O}(\lambda^2)\mu^2\varphi^2+\mathcal{O}_\mu(\lambda)\varphi^2\]\theta^2 |Y|^2dxdt.
\end{eqnarray*}
Also, notice that
\begin{eqnarray*}
&&2\Phi zdz=4\Delta\ell zdz=2\Delta \ell[d(z^2)-(dz)^2]\\
&&=2d(\Delta \ell z^2)-2\Delta\ell_t z^2dt-2\Delta\ell(dz)^2.
\end{eqnarray*}
This implies that
\begin{equation}\label{fl10}
\dbE\int_Q 2\Phi zdzdx=
\dbE\int_Q\Big[\mathcal{O}(\lambda)\mu^2\varphi^2|z|^2
+\mathcal{O}(\lambda)\mu^2\varphi\theta^2|Y|^2\Big]dxdt.
\end{equation}

Combining (\ref{fl5})-(\ref{fl10}) with (\ref{fl00}),  one can get the desired inequality (\ref{*}). \endpf


\end{document}